\documentclass[10pt]{amsart}
 
\usepackage{amscd,a4,graphics,amssymb}
\usepackage[all]{xy}
\usepackage{color}

\newtheorem{thm}{Theorem}[section] \newtheorem{lemma}[thm]{Lemma}
\newtheorem{prop}[thm]{Proposition}
\newtheorem{cor}[thm]{Corollary}

\theoremstyle{definition}
 \newtheorem{dfn}[thm]{Definition}
 \newtheorem{rmk}[thm]{Remark}
\newenvironment{pf}{\medskip\noindent{{\em Proof: }}}{\qed}

\newcommand{\Spec}{\mathrm{Spec}}
\newcommand{\Frac}{\mathrm{Frac}}
 
\newcommand{\Char}{\mathrm{char}}
\newcommand{\Sep}{\mathrm{sep}}
\newcommand{\Red}{\mathrm{red}}
\newcommand{\SNC}{\mathrm{SNC}}
\newcommand{\Gal}{\mathrm{Gal}}
\newcommand{\Md}{\mathrm{md}}
\newcommand{\Gcd}{\mathrm{gcd}}
\newcommand{\Deg}{\mathrm{deg}}
\newcommand{\Mult}{\mathrm{mult}}
\newcommand{\Lcm}{\mathrm{lcm}}
\newcommand{\Min}{\mathrm{min}}
\newcommand{\Ord}{\mathrm{ord}}
\newcommand{\Sing}{\mathrm{sing}}
\newcommand{\Id}{\mathrm{id}}

\newcommand{\Mod}{\mathrm{mod}}

\begin{document}


  \title{Stable reduction of curves and tame ramification}


  \author{Lars Halvard Halle}\email{halle@math.kth.se}
  \address{Department of Mathematics, KTH, S--100 44 Stockholm,
    Sweden}


\begin{abstract}
We study stable reduction of curves in the case where a tamely ramified base extension is sufficient. If $X$ is a smooth curve defined over the fraction field of a strictly henselian discrete valuation ring, there is a criterion, due to T. Saito, that describes precisely, in terms of the geometry of the minimal model with strict normal crossings of $X$, when a tamely ramified extension suffices in order for $X$ to obtain stable reduction. For such curves we construct an explicit extension that realizes the stable reduction, and we furthermore show that this extension is minimal. We also obtain a new proof of Saito's criterion, avoiding the use of $\ell$-adic cohomology and vanishing cycles.
\end{abstract}


\thanks{Work partly done during the Moduli Spaces-year 2006/2007 at Institut Mittag-Leffler (Djursholm, Sweden)}

\maketitle

\section{Introduction}

Let us first recall the so-called \emph{stable reduction} theorem.

\begin{thm}[Deligne-Mumford, \cite{DelMum}]
Let $ R $ be a discrete valuation ring, and let $X$ be a smooth, projective, geometrically connected curve of genus $ g(X) \geq 2 $ over $K=\Frac(R)$. Then there exists a finite, separable extension $ K \subset L $ such that $X \otimes_K L$ has stable reduction over $R_L$, the integral closure of $R$ in $L$. That is, there exists a stable curve $ \mathcal{X}_L $ over $\Spec(R_L)$ such that the generic fiber is isomorphic to $X \otimes_K L$.
\end{thm}

If the residue characteristic is zero, one can find explicit extensions of $R$ that realize stable reduction for $X$. This is done by considering a suitable regular model with normal crossings $ \mathcal{Y} $ for $X$ over $R$, and taking an extension of discrete valuation rings $R'/R$ of ramification index divisible by the multiplicities of all irreducible components in the special fiber of $ \mathcal{X} $. One can then show that the normalization of the pullback $ \mathcal{X} \times_{\Spec(R)} \Spec(R') $ has only $ A_n $-singularitites, which can be explicitly resolved, and that the minimal desingularization is semi-stable. The stable model is then obtained by contracting all $(-2)$-curves in the special fiber (cf.~\cite{HarMor}, Proposition 3.39, or for a  more general result \cite{Liubook}, Proposition 10.4.6). Furthermore, having this description is often very useful when one wants to compute the \emph{stable reduction} of $X$, that is, the special fiber of the stable model.

In positive characteristic, it can often be hard to find explicit extensions that realize stable reduction. One of the purposes of this paper is, in those cases where a \emph{tamely ramified} extension suffices, to show that the geometry of the special fiber of a suitable normal crossings model $ \mathcal{X} $ for $X$ over $R$ still contains enough information so that we can find explicit extensions over which $X$ obtains stable reduction. To do this, we show that the stable model can be constructed using a certain base-change/normalization/desingularization/contraction procedure, generalizing the one above. The main problem that needs to be overcome is the fact that components in the special fiber of $ \mathcal{X} $ may have multiplicities divisible by $p$.

In \cite{Saito}, T. Saito gave a new proof of the stable reduction theorem, using $\ell$-adic cohomology and the theory of vanishing cycles, relating the monodromy action of the Galois group $ \Gal(K_{\Sep}/K) $ on $ H^1(X_{K_{\Sep}}, \mathbb{Q}_{\ell}) $ with the geometry of certain normal crossings models for $X$ over $R$. Furthermore, he also gave a geometric criterion for when the \emph{wild ramification group} $ P \subset \Gal(K_{\Sep}/K) $ acts trivially on $ H^1(X_{K_{\Sep}}, \mathbb{Q}_{\ell}) $ (\cite{Saito}, Theorem 3). It may be seen that $P$ acts trivially in the sense above if and only if $X$ has stable reduction after a tamely ramified base extension. Having this geometric description will be crucial for this paper. A precise formulation is given in Theorem \ref{S-crit} below. In particular, Theorem \ref{S-crit} describes precisely the components that may have multiplicity divisible by the residue characteristic. 

If $R$ is strictly henselian, with algebraically closed residue field, it is known that there exists a finite extension $ K \subset L \subset K_{\Sep} $ of $K$ minimal with the property that $ X_L = X \otimes_{K} L $ has stable reduction over $R_L$, where $R_L$ is the integral closure of $R$ in $L$ (\cite{Liubook}, Theorem 10.4.44). In the case where $X$ has stable reduction after a tame extension, we determine exactly this minimal extension $L/K$, and thus we generalize a result by G. Xiao in characteristic zero (\cite{Xiao}, Proposition 1).

\subsection{Notation} We list here some notation that will be valid throughout the text.

$R = $ a discrete valuation ring, with uniformizing parameter $\pi$.

$k = $ the residue field of $R$, assumed to be algebraically closed.

$p = \Char(k) $.

$K = $ the fraction field of $R$.

$K_{\Sep} = $ the separable closure of $K$. 

$S = \Spec(R)$.

$ X = $ a smooth, projective and geometrically connected curve over $K$, with genus $ g(X) \geq 2 $.

\subsection{Strict normal crossings models of $X$}
It is well-known that we can extend $X$ to a relative curve $ \mathcal{X} $ over $S$ in such a way that $ \mathcal{X} $ is a regular surface, the irreducible components of the special fiber $ \mathcal{X}_k $ are smooth and such that $ (\mathcal{X}_k)_{\Red} $ is a \emph{strict} normal crossings divisor. Such a surface will be called an $ \SNC $-model for $X$. Furthermore, we can choose $ \mathcal{X} $ minimal with respect to these properties (cf.~\cite{Liubook}, Prop. 9.3.36).

\subsection{Saito's criterion}

\begin{thm}[\cite{Liubook}, Theorem 10.4.47] \label{S-crit}
Let $X/K$ be as above, and let $\mathcal{X}$ be the minimal $\SNC$-model of $X$ over $R$, where $R$ is strictly henselian. The following conditions are equivalent:
\begin{enumerate}
\item The minimal extension $K \subset L$ that realizes the stable reduction of $ X $ is tamely ramified.

\item Every irreducible component $C$ of $\mathcal{X}_k$ whose multiplicity in $\mathcal{X}_k$ is divisible by $p$ satisfies the following condition $(*)$ :

$(*)$ $C$ is isomorphic to $\mathbb{P}^1_k$, and intersects the other components of $\mathcal{X}_k$ in exactly two points and the components that meet $C$ have multiplicities in $\mathcal{X}_k$ that are not divisible by $p$.
\end{enumerate}
\end{thm}

We shall refer to statement (2) in Theorem \ref{S-crit} as \emph{Saito's criterion}. Furthermore, if $\mathcal{X}/S$ is an SNC-model for $X$, and the special fiber $\mathcal{X}_k$ satisfies Saito's criterion, we shall say that $\mathcal{X}$ satisfies Saito's criterion.

\begin{rmk}
Recall that when $R$ is henselian, and $n$ is an integer not divisible by $p$, then the extension $ K \rightarrow K' := K[\pi']/(\pi'^n - \pi) $, where $\pi$ is a uniformizing parameter of $R$, is tamely ramified over $(\pi)$ only, and the integral closure $R'$ of $R$ in $K'$ is $ R[\pi']/(\pi'^n - \pi) $, which is a discrete valuation ring with uniformizing parameter $\pi'$. Conversely, if $ K \rightarrow K' $ is a tamely ramified extension of degree $n$, ramified only over $(\pi)$, we can write it in the form $ K' = K[\pi']/(\pi'^n - \pi) $, and the integral closure of $R$ in $K'$ is $ R[\pi']/(\pi'^n - \pi) $ (cf. ~\cite{Neukirch}, Proposition II.7.7). In this paper, whenever we say that $ \Spec(R') \rightarrow \Spec(R) $ is a tamely ramified extension, we shall mean that it is of this type.
\end{rmk}

\subsection{Overview}
We give here a short overview of this paper.

$\bullet$ In Section 2, we consider the following situation: Let $\mathcal{X}/S$ be an SNC-model, where $S$ is the spectrum of a \emph{complete} discrete valuation ring, and let $ \mathcal{X}' $ be the normalization of $ \mathcal{X} \times_S S' $, where $ S'/S $ is a tamely ramified extension. Under certain assumptions on the special fiber of $\mathcal{X}$, we compute the completion of the local ring of any closed point in the special fiber of $ \mathcal{X}'/S' $.

$\bullet$ In Section 3, we study the induced morphism $ \mathcal{X}' \rightarrow \mathcal{X} $, and the special fiber of ~$ \mathcal{X}' $.

$\bullet$ In Section 4, we introduce \emph{tame cyclic quotient singularities}, following \cite{CED}, and show that $ \mathcal{X}' $ has at worst such singularities. Furthermore, we describe the minimal desingularization $ \mathcal{X}'_{\Md} \rightarrow \mathcal{X}' $.

$\bullet$ Sections 5 and 6 consist of a combinatorial study of the special fiber of $ \mathcal{X}'_{\Md} $, with emphasis on the contraction of smooth and rational components.

$\bullet$ Finally, in Section 7, in the case where $\mathcal{X}/S$ is minimal, and satisfies Saito's criterion, we find an explicit tamely ramified extension that realizes stable reduction for $X/K$. This extension depends only on the geometry of $\mathcal{X}/S$. Furthermore, we show that this is the minimal extension realizing stable reduction. As a corollary, we obtain a new and more geometric proof of Theorem \ref{S-crit}, without the use of vanishing cycles. 

\subsection{Acknowledgments}
I would like to thank Dino Lorenzini for useful suggestions. I would also like to thank my thesis advisor Carel Faber for discussing the material in this paper with me.

\section{Local computations}\label{section 2}
\subsection{Setup}\label{section 2.1}
We will, unless otherwise mentioned, assume throughout the rest of the paper that $R$ is complete. Let $n$ be an integer not divisible by $p$. Let $ K' = K[\pi']/(\pi'^n - \pi) $, and let $ R' $ be the integral closure of $R$ in $K'$. Then $ R' $ is a complete discrete valuation ring, finite over $R$ of ramification index $n$. Let $ S' = \Spec(R') $, and consider the diagram
$$ \xymatrix{
 \mathcal{X}' \ar[d] \ar[r]  & \mathcal{X}_{S'} \ar[d] \ar[r] &  \mathcal{X} \ar[d] \\
 S' \ar[r]^{\Id} & S' \ar[r] & S. } $$
The pullback $ \mathcal{X}_{S'} := \mathcal{X} \times_S S' $ is flat over $S'$, with smooth and irreducible generic fiber, hence it is integral (\cite{Liubook}, Prop. 4.3.8). Furthermore, $ \mathcal{X}_{S'} $ is excellent, since it is of finite type over the excellent scheme $S'$, so the normalization $ \mathcal{X}' \rightarrow \mathcal{X}_{S'} $ is finite (\cite{Liubook}, Theorem 8.2.39). Therefore, the composition $ f : \mathcal{X}' \rightarrow \mathcal{X} $ is a finite morphism.

\subsection{Local rings, completion and normalization}
Let $ x \in \mathcal{X}_k $ be a closed point, and let $ \mathcal{O}_{\mathcal{X},x} $ be the local ring of $ \mathcal{X} $ at $x$. The ring $ \mathcal{O}_{\mathcal{X},x} \otimes_R R' $ is the local ring of $ \mathcal{X}_{S'} $ at the unique point mapping to $x$, and hence is reduced and excellent. Let $ (\mathcal{O}_{\mathcal{X},x} \otimes_R R')' $ denote the normalization of $ \mathcal{O}_{\mathcal{X},x} \otimes_R R' $ in its total ring of fractions. Then we have that $ (\mathcal{O}_{\mathcal{X},x} \otimes_R R')' $ is semi-local, and the maximal ideals correspond to the points $ x'_1, \ldots , x'_{m} $ of $ \mathcal{X}' $ mapping to $x$ via $f$. The localization in a maximal ideal is the local ring $ \mathcal{O}_{\mathcal{X}',x'_i} $ of $ \mathcal{X}' $ at $ x'_i $, for some $ i \in \{1, \ldots, m\} $. The induced homomorphism $ \mathcal{O}_{\mathcal{X},x} \rightarrow \mathcal{O}_{\mathcal{X}',x'_i} $ may be identified with the local homomorphism induced by $f$.

Since $ \mathcal{O}_{\mathcal{X},x} \rightarrow \mathcal{O}_{\mathcal{X},x} \otimes_R R' $ is finite, tensoring with the completion $ \mathcal{O}_{\mathcal{X},x} \rightarrow \widehat{\mathcal{O}}_{\mathcal{X},x} $ gives a cartesian diagram
$$ \xymatrix{
\mathcal{O}_{\mathcal{X},x} \ar[d] \ar[r] & \mathcal{O}_{\mathcal{X},x} \otimes_R R' \ar[d]\\
\widehat{\mathcal{O}}_{\mathcal{X},x} \ar[r] & \widehat{\mathcal{O}}_{\mathcal{X},x} \otimes_R R',}
$$
where $ \widehat{\mathcal{O}}_{\mathcal{X},x} \otimes_R R' $ is the completion of $ \mathcal{O}_{\mathcal{X},x} \otimes_R R' $, and $ \widehat{\mathcal{O}}_{\mathcal{X},x} \rightarrow \widehat{\mathcal{O}}_{\mathcal{X},x} \otimes_R R' $ is the completion of the homomorphism $ \mathcal{O}_{\mathcal{X},x} \rightarrow \mathcal{O}_{\mathcal{X},x} \otimes_R R' $.

The ring $ \widehat{\mathcal{O}}_{\mathcal{X},x} \otimes_R R' $ is reduced (\cite{Abbes}, Lemme A.4). Let $ (\widehat{\mathcal{O}}_{\mathcal{X},x} \otimes_R R')' $ denote the normalization in its total ring of fractions. We have that
$$ (\widehat{\mathcal{O}}_{\mathcal{X},x} \otimes_R R')' \cong C((\mathcal{O}_{\mathcal{X},x} \otimes_R R')'), $$
where $ C((\mathcal{O}_{\mathcal{X},x} \otimes_R R')') $ denotes the completion with respect to the radical (\cite{EGAIV}, 7.8.3 (vii)). On the other hand, we also have that 
$$ C((\mathcal{O}_{\mathcal{X},x} \otimes_R R')') \cong \prod_{i = 1}^{m} \widehat{\mathcal{O}}_{\mathcal{X}',x'_i}. $$ 
Therefore, the compositions
$$ \widehat{\mathcal{O}}_{\mathcal{X},x} \rightarrow \widehat{\mathcal{O}}_{\mathcal{X},x} \otimes_R R' \rightarrow ( \widehat{\mathcal{O}}_{\mathcal{X},x} \otimes_R R')' \cong \prod_{i=1}^m \widehat{\mathcal{O}}_{\mathcal{X}',x'_i} \rightarrow \widehat{\mathcal{O}}_{\mathcal{X}',x'_i}, $$
where the last map is the projection onto the $i$-th factor, describe the maps of the completed local rings induced by $ \mathcal{X}' \rightarrow \mathcal{X} $. In what follows, we shall make these maps more explicit.

\subsection{The local rings of $ \mathcal{X}' $}
Let $ x' \in \mathcal{X}'_{k} $ be a closed point mapping to $ f(x') = x \in \mathcal{X}_k $. We will make the assumption that either
\begin{enumerate}
\item $x$ belongs to a unique component of $ \mathcal{X}_k $, or

\item $x$ is an intersection point of two distinct components of $ \mathcal{X}_k $, where at least one of the components has multiplicity not divisible by $p$.
\end{enumerate}
Under this assumption, we shall in the following compute $ \widehat{\mathcal{O}}_{\mathcal{X}',x'} $. We will treat the two cases above independently. The local analytic structure of $ \mathcal{X}' $ at $x'$ will only depend on the structure of $ \mathcal{X} $ at $x$.

\subsection{One branch}\label{section 2.7}
In case (1) above, $x$ belongs to a unique component of $ \mathcal{X}_k $, and we can find an isomorphism
$$ \widehat{\mathcal{O}}_{\mathcal{X},x} \cong R[[u,v]]/(\pi - c_0 v^b), $$ for some unit $ c_0 \in R[[u,v]] $ (cf. ~\cite{CED}, proof of Lemma 2.3.2). Let $ b=b'l $ and $ n = n'l $, where $ l = \Gcd(b,n) $. Since $n$ is not divisible by $p$, we can find a unit $c \in R[[u,v]]$ such that $ c^{ln'} = c_0 $. Then we have 

$$ \{ \widehat{\mathcal{O}}_{\mathcal{X},x} \otimes_R R' \}' \cong \prod_{\eta \in \boldsymbol{\mu}_l} \{ R'[[u,v]]/(\pi'^{n'} - \eta c^{n'} v^{b'}) \}'. $$

The factors $ R'[[u,v]]/(\pi'^{n'} - \eta c^{n'} v^{b'}) $ are reduced, since $ \widehat{\mathcal{O}}_{\mathcal{X},x} \otimes_R R' $ is reduced. After possibly taking an $n'$-th root of $ \eta$, it suffices to compute the normalization of the ring $ R'[[u,v]]/(\pi'^{n'} - c^{n'} v^{b'}) $. Consider the $R'$-homomorphism
$$ \Phi : R'[[u,v]] \rightarrow R'[[s,t]], $$
defined by $ (u,v) \mapsto (s,t^{n'}) $. Then we have that 
$$ \pi'^{n'} - c^{n'} v^{b'} = \prod_{ \xi \in \boldsymbol{\mu}_{n'} } (\pi' - \xi c t^{b'}) \in R'[[s,t]]. $$
Each of the factors $ R'[[s,t]]/(\pi' - \xi c t^{b'}) $ is regular, and hence irreducible. Furthermore, $ \boldsymbol{\mu}_{n'} $ acts on $ R'[[s,t]] $ by $ [ \sigma ](t) = \sigma t $, for any $ \sigma \in \boldsymbol{\mu}_{n'} $, and the invariant ring for this action is $R'[[u,v]]$. The factors $ \pi' - \xi c t^{b'} $ are permuted under this action, and since $ \Gcd(n',b') = 1 $, it is easily seen that there is only one orbit. Consequently, $ V(\pi'^{n'} - c^{n'} v^{b'}) \subset \Spec(R'[[u,v]]) $ is irreducible as well as reduced. So the ideal $ (\pi'^{n'} - c^{n'} v^{b'}) \subset R'[[u,v]] $ is a prime ideal.

The homomorphism $ \Phi $ induces an injective $R'$-homomorphism
$$ \phi : A := R'[[u,v]]/(\pi'^{n'} -  c^{n'} v^{b'}) \rightarrow B := R'[[s,t]]/(\pi' - c t^{b'}), $$
and the elements $ 1, t, \ldots , t^{n'-1} $ generate $B$ as an $A$-module. Since $ \Gcd(b',n') = 1 $, we can find integers $ \alpha, \beta $, such that  $ \alpha b' + \beta n' = 1 $. Furthermore, we have the relations $ t^{n'} = v $ and $ t^{b'} = c^{-1} \pi' $ in $B$, so $ t^{n'} $ and $ t^{b'} $ lie in the image of $ \phi $. But this implies that 
$$ (t^{b'})^{\alpha}(t^{n'})^{\beta} = t \in \Frac(A), $$
so it follows that $ \Frac(A) \cong \Frac(B) $, and therefore $ A' \cong B $.

We sum this up in the following proposition.

\begin{prop}\label{prop. 2.8}
Let $ x \in \mathcal{X}_k $ be a closed point with $ \widehat{\mathcal{O}}_{\mathcal{X},x} \cong R[[u,v]]/(\pi- c_0 v^b) $. If $ x' \in \mathcal{X}'_{k} $ is a closed point mapping to $x$, then we have that 
$$ \widehat{\mathcal{O}}_{\mathcal{X}', x'} \cong R'[[s,t]]/(\pi' - c t^{b'}), $$
where $ b' = b/\Gcd(b,n) $.
\end{prop}

\begin{cor}\label{cor. 2.9}
$ \mathcal{X}' $ is regular at $x'$. Furthermore, there is exactly one irreducible component $ C' $ of $ \mathcal{X}'_{k} $ passing through $x'$, and this branch is smooth at $x'$.
\end{cor}
\begin{pf}
By Proposition \ref{prop. 2.8}, $ \widehat{\mathcal{O}}_{\mathcal{X}', x'} $ is regular, and therefore also $ \mathcal{O}_{\mathcal{X}', x'} $ is regular.

Since the completion  $ \mathcal{O}_{\mathcal{X}', x'} \rightarrow \widehat{\mathcal{O}}_{\mathcal{X}', x'} $ is faithfully flat, it follows that there is exactly one irreducible component $C'$ of the special fiber passing through $x'$. Let $ I $ be the ideal of $ C' $ in $ \mathcal{O}_{\mathcal{X}', x'} $. As $ I $ is a prime ideal, we have that $ \mathcal{O}_{C', x'} = \mathcal{O}_{\mathcal{X}', x'}/I $ is an integral domain, and in particular reduced. Furthermore, since $ \mathcal{O}_{C', x'} $ is excellent, the completion
\begin{equation}\label{equation 1} 
\widehat{\mathcal{O}}_{C', x'} \cong \widehat{(\mathcal{O}_{\mathcal{X}', x'}/I)} \cong \widehat{\mathcal{O}}_{\mathcal{X}', x'}/ I \cdot \widehat{\mathcal{O}}_{\mathcal{X}', x'} 
\end{equation}
is also reduced (\cite{Liubook}, Proposition 8.2.41).

As $ V(I) \subset \Spec(\widehat{\mathcal{O}}_{\mathcal{X}', x'}) $ is obviously irreducible, it follows that $ I \cdot \widehat{\mathcal{O}}_{\mathcal{X}', x'} $ is a prime ideal, and then we necessarily get that
$$ I \cdot \widehat{\mathcal{O}}_{\mathcal{X}', x'} = (t) \subset \widehat{\mathcal{O}}_{\mathcal{X}', x'}. $$

From Equation \ref{equation 1} above, it then follows that
$$ \widehat{\mathcal{O}}_{C', x'} \cong k[[s]], $$
so $C'$ is indeed smooth at $x'$.
\end{pf}

\subsection{Two branches}
We consider now case (2), where $x$ is an intersection point of two distinct components of $ \mathcal{X}_k $. Then we can find an isomorphism
$$ \widehat{\mathcal{O}}_{\mathcal{X},x} \cong R[[u,v]]/(\pi-u^av^b), $$ 
where $a$ and $b$ are the multiplicities of the components meeting at $x$ (cf. ~\cite{CED}, proof of Lemma 2.3.2). Here the assumption that $ p $ does not divide both $a$ and $b$ is necessary to get this easy \emph{polynomial} form. This will be important when we consider the desingularization of $ \mathcal{X}' $ at $x'$.

Let us write $ a = a'd $, $ b= b'd $ and $ n = n'd $, where $ d = \Gcd(a,b,n) $. Then we have that
$$ (\widehat{\mathcal{O}}_{\mathcal{X},x} \otimes_R R')' \cong \prod_{\xi \in \boldsymbol{\mu}_d} \{ R'[[u,v]]/(\pi'^{n'} - \xi u^{a'}v^{b'}) \}'. $$
Let $ e = \Gcd(n',a') $ and $ c = \Gcd(n',b') $. Note that in particular $ \Gcd(e,c) = 1 $. We can now write $ n' = n''ce $, $ b' = cb'' $ and $ a' = ea'' $.

If $ x' \in \mathcal{X}'_{k} $ is a closed point mapping to $ x $, then we have that 
$$ \widehat{\mathcal{O}}_{\mathcal{X}',x'} \cong \{R'[[u,v]]/(\pi'^{n'} - \xi u^{a'}v^{b'})\}', $$
for some $ \xi \in \boldsymbol{\mu}_d $. After a change of coordinates, we may assume that $ \xi = 1 $. In order to normalize $ R'[[u,v]]/(\pi'^{n'} - \xi u^{a'}v^{b'}) $, we first consider the $R'$-algebra homomorphism
$$ \Psi : R'[[u,v]] \rightarrow R'[[s,t]] $$
given by $ (u,v) \mapsto (s^c,t^e) $. We have that
$$ \pi'^{n'} - u^{a'}v^{b'} = \prod_{\xi \in \boldsymbol{\mu}_{e c}} (\pi'^{n''} - \xi s^{a''} t^{b''}) \in R'[[s,t]]. $$

\begin{lemma}\label{prime}
The element $ \pi'^{n''} - s^{a''} t^{b''} \in R'[[s,t]]$ generates a prime ideal. 
\end{lemma}
\begin{pf}
Let the $R'$-homomorphism $ R'[[s,t]] \rightarrow R'[[z,w]] $ be defined by~$ s \mapsto z^{n''} $ and $ t \mapsto w^{n''} $. We have that 
$$ \pi'^{n''} - s^{a''} t^{b''} = \prod_{\xi \in \boldsymbol{\mu}_{n''}} (\pi' - \xi z^{a''} w^{b''}) \in R'[[z,w]]. $$
It is easily seen that $ V(\prod_{\xi \in \boldsymbol{\mu}_{n''}} (\pi' - \xi z^{a''} w^{b''})) \subset \Spec(R'[[z,w]]) $ is reduced, and so it follows that $ V(\pi'^{n''} - s^{a''} t^{b''}) \subset \Spec(R'[[s,t]]) $ is reduced. Furthermore, let $ \boldsymbol{\mu}_{n''} \times \boldsymbol{\mu}_{n''} $ act on $ R'[[z,w]] $ by $ (\xi_1,\xi_2)[z] = \xi_1 z $ and $ (\xi_1,\xi_2)[w] = \xi_2 w $. The invariant ring for this action is $ R'[[s,t]] $. The schemes $ V(\pi' - \xi z^{a''} w^{b''}) $ are regular, and hence irreducible, and are easily seen to belong to the same orbit under this action. It follows that $ V(\pi'^{n''} - s^{a''} t^{b''}) \subset \Spec(R'[[s,t]]) $ is irreducible. Therefore, we have that $ (\pi'^{n''} - s^{a''} t^{b''}) \subset R'[[s,t]] $ is a prime ideal.
\end{pf}

\vspace{0.5cm}

In a similar way as in Lemma \ref{prime}, we can let $ \boldsymbol{\mu}_{e c} = \boldsymbol{\mu}_{e} \times \boldsymbol{\mu}_{c} $ act on $ R'[[s,t]] $, and show that $ \pi'^{n'} - u^{a'}v^{b'} \in R'[[u,v]] $ generates a prime ideal. It follows that $ \Psi $ induces an injective homomorphism
$$ A := R'[[u,v]]/(\pi'^{n'} - u^{a'}v^{b'}) \rightarrow B := R'[[s,t]]/(\pi'^{n''} - s^{a''}t^{b''}), $$
where $ (u,v) \mapsto (s^c,t^e) $. Furthermore, $B$ is a finite $A$-module, generated by the finitely many elements $ s^it^j $, where $ 0 \leq i < c $, $ 0 \leq j < e $.

\begin{lemma}\label{lemma 2.1}
We have that $ \Frac(A) = \Frac(B) $, and hence the normalization of $A$ equals the normalization of $B$. (So $B$ is a partial normalization of $A$).
\end{lemma}
\begin{pf}
The elements $ s^c $, $ t^e $ and $ s^{a''} t^{b''} $ lie in the image of $A$. Since $ \Gcd(c,e)=1 $, there exist integers $ \alpha $ and $ \beta $ such that $ \alpha c + \beta e = 1 $. But then we get that
$$ (s^{a''}t^{b''})^{\beta e } (s^c)^{a'' \alpha} = (s^{a''})^{\alpha c + \beta e} (t^e)^{\beta b''} = s^{a''} (t^e)^{\beta b''} \in \Frac(A). $$ 
It follows that also $ s^{a''} \in \Frac(A) $. But $ \Gcd(a'',c)=1 $, so there exist integers $ \alpha_1 $ and $ \beta_1 $ such that $ \alpha_1 a'' + \beta_1 c = 1 $. Consequently 
$$ s = (s^{a''})^{\alpha_1} (s^c)^{\beta_1} \in \Frac(A). $$

Arguing in a similar way, we find that also $ t \in \Frac(A) $. Hence $B$ is generated as an $A$-module by finitely many elements that lie in $\Frac(A)$, and so the result follows.
\end{pf}

\subsection{Group action}\label{section 2.5}
Consider the ring $ C := R'[[z,w]]/(\pi'-z^{a''}w^{b''}) $. From the proof of Lemma \ref{prime} it follows that the $R'$-homomorphism
$$ B = R'[[s,t]]/(\pi'^{n''} - s^{a''}t^{b''}) \rightarrow C = R'[[z,w]]/(\pi'-z^{a''}w^{b''}), $$
where $ (s,t) \mapsto (z^{n''},w^{n''}) $, is injective. Furthermore, $C$ is a finite $B$-module, generated by the elements $z^iw^j$, where $ 0 \leq i < n'' $, $ 0 \leq j < n'' $.

Since $ a'' $ and $b''$ are relatively prime to $n''$, we can find a unit $ r \in (\mathbb{Z}/n'')^* $ such that $ r b'' + a'' \equiv 0 $ mod $n''$. Then we let $ G = \boldsymbol{\mu}_{n''} $ act on $C$ by $ [\xi](z) = \xi z $ and $ [\xi](w) = \xi^r w $, for any $ \xi \in \boldsymbol{\mu}_{n''} $.

Since $C$ is regular, the invariant ring under the $ G $-action is a normal complete local ring, and we shall see that this is indeed the normalization of $B$, and hence of $A$. We first prove that $ C^{G} $ is a finite $B$-module, and find an explicit set of generators.

\begin{lemma}\label{lemma 2.3}
The invariant ring of $C$ under the action of $ G $ is generated as a $B$-module by the $G$-invariant monomials of the form $ z^iw^j $, where $ 0 \leq i,j < n'' $.
\end{lemma}
\begin{pf}
Let $\xi$ be a primitive $n''$-th root of unity in $R'$. If $ z^iw^j $ is invariant under $G$, we obviously have that $ [\xi^k](z^iw^j) = z^iw^j $ for all $ 0 \leq k < n'' $, and consequently
$$ z^iw^j + [\xi](z^iw^j) + \ldots + [\xi^k](z^iw^j) + \ldots + [\xi^{n''-1}](z^iw^j) = n'' z^iw^j. $$

Assume now that $ z^i w^j $ is \emph{not} invariant under the $G$-action. Then we definitely have that $ i + r j = n''N + r' $, where $ 0 < r' < n'' $. Furthermore, if $ k $ is an integer such that $ 0 \leq k < n'' $, then
$$ [\xi^k](z^i w^j) = (\xi^k z)^i ((\xi^k)^r w)^j = \xi^{k(i + r j)} z^i w^j = \xi^{kr'} z^i w^j. $$

So it follows that 
$$ z^iw^j + [\xi](z^iw^j) + \ldots + [\xi^k](z^iw^j) + \ldots + [\xi^{n''-1}](z^iw^j) $$
$$ = (1 + \xi^{r'} + \ldots + \xi^{kr'} + \ldots + \xi^{(n''-1)r'} ) z^i w^j.  $$

But $ 1 + \xi^{r'} + \ldots + \xi^{kr'} + \ldots + \xi^{(n''-1)r'} = 0 $, since
 $$ (1-\xi^{r'}) (1 + \xi^{r'} + \ldots + \xi^{kr'} + \ldots + \xi^{(n''-1)r'} ) = 0 $$ 
in $R'$, which is an integral domain, and $ (1-\xi^{r'}) \neq 0 $, since $ \xi $ is a primitive root and $ 0 < r' < n'' $.

If now $ F = \sum_{0 \leq i,j < n''} f_{i,j} z^iw^j $, where $ f_{i,j} \in B $, is an element in $C$ which is invariant under $G$, we have that
$$ \sum_{k=0}^{n''-1} [\xi^k](F) = F + [\xi](F) + \ldots + [\xi^k](F) + \ldots + [\xi^{n''-1}](F) = n'' F. $$
On the other hand, by the computations above, we have that
$$ \sum_{k=0}^{n''-1} [\xi^k](F) = \sum_{0 \leq i,j < n''} f_{i,j} z^iw^j + \ldots + [\xi^{n''-1}] (\sum_{0 \leq i,j < n''} f_{i,j} z^iw^j) = n'' \sum_{i',j'} f_{i',j'} z^{i'}w^{j'}, $$
where the last sum runs over those $ 0 \leq i',j' < n'' $ such that  $ z^{i'}y^{j'} $ is \emph{invariant} under the action of $ G $.
\end{pf}

\begin{lemma}\label{lemma 2.2}
If a monomial $ z^{e_1}w^{e_2} \in C $ is invariant under the $ G $-action, then $ z^{e_1}w^{e_2} \in \Frac(B) $. 
\end{lemma}
\begin{pf}
Let us first make the observation that if $ z^{e_1}w^{e_2} $ is invariant for the action of $G$, then we have that $ z^{e_1}w^{e_2} = [\xi](z^{e_1}w^{e_2}) = z^{e_1}w^{e_2} \xi^{e_1+r e_2} $, for any $ \xi \in G $. Hence $ e_1+r e_2 = n'' N $ for some integer $N$. By construction we have that $ rb'' = n''M - a'' $ for some integer $M$. Note also that the elements $ z^{a''}w^{b''} $, $ z^{n''} $ and $w^{n''}$ all lie in the image of $B$. But now we compute
$$ (z^{a''}w^{b''})^{e_1} = z^{a'' e_1} w^{b''n''N-rb''e_2} = z^{a'' e_1} w^{b''n''N - e_2(n''M-a'')} = $$
$$ z^{a'' e_1} w^{a''e_2} (w^{n''})^{b''N-e_2M} = (z^{e_1}w^{e_2})^{a''} (w^{n''})^{b''N-e_2M}, $$
which implies that $ (z^{e_1}w^{e_2})^{a''} \in \Frac(B) $. On the other hand, we have that $ (z^{e_1}w^{e_2})^{n''} \in \Frac(B) $. Since $ \Gcd(a'',n'') = 1 $, we can find integers $ \alpha, \beta $ such that $ \alpha a'' + \beta n'' = 1 $. From this, it follows that 
$$ (z^{e_1}w^{e_2})^{ \alpha a''} (z^{e_1}w^{e_2})^{\beta n''} = z^{e_1}w^{e_2}, $$
which finally gives that $ z^{e_1}w^{e_2} \in \Frac(B) $. 
\end{pf}

\subsection{Normalization}\label{section 2.6}

\begin{prop}
The normalization of $A$ is the invariant ring $C^G$ of $C$ under the action of $G = \boldsymbol{\mu}_{n''} $ introduced in Section \ref{section 2.5}.
\end{prop}
\begin{pf}
By Lemma \ref{lemma 2.3}, we know that $C^G$ is finite over $B$, generated by finitely many monomials that are invariant under $G$. On the other hand, Lemma \ref{lemma 2.2} shows that each of these monomials lies in $\Frac(B)$, so it follows that $ \Frac(C^G) = \Frac(B) $. But $C^G$ is \emph{normal}, and must therefore equal the normalization of $B$. Lemma \ref{lemma 2.1} then shows that we may identify $ A' \cong C^G $.
\end{pf}

\vspace{0.5cm}

We sum up the results regarding the normalization in the proposition below.

\begin{prop}\label{prop. 2.5}
Let $ x \in \mathcal{X}_k $ be a closed point with $ \widehat{\mathcal{O}}_{\mathcal{X},x} \cong R[[u,v]]/(\pi-u^av^b) $. If $ x' \in \mathcal{X}'_{k} $ is a closed point mapping to $x$, then we have that 
$$ \widehat{\mathcal{O}}_{\mathcal{X}',x'} \cong \{ R'[[z,w]]/(\pi'-z^{a''}w^{b''}) \}^{\boldsymbol{\mu}_{n''}}, $$
where $ a'' = a/\Gcd(a,n) $, $ b'' = b/\Gcd(b,n) $, $ n'' = (n \cdot \Gcd(a,b,n))/(\Gcd(a,n) \cdot \Gcd(b,n)) $, and where the $ \boldsymbol{\mu}_{n''} $-action is given as in Section \ref{section 2.5}. 
\end{prop}

We can now describe the irreducible components of the special fiber of $ \mathcal{X}' $ locally at $x'$.
 \begin{prop}\label{prop. 2.6}
Let us keep the hypotheses from Proposition \ref{prop. 2.5}. There are precisely two irreducible components of $ \mathcal{X}'_k $ passing through $x'$. Furthermore, these components are smooth at $x'$.
\end{prop}
\begin{pf}
Consider first $ C = R'[[z,w]]/(\pi'-z^{a''}w^{b''}) $, with the $G$-action given in Section \ref{section 2.5}. The special fiber of $C$ has two irreducible components, with ideals $(z)$ and $(w)$. As these are stable under the $G$-action, it follows there are two distinct irreducible components of the special fiber of $C^G$, namely the images of $(z)$ and $(w)$.

Let $ I \subset C $ be any ideal that is stable under $G$. Since the order of $G$ is invertible in $C$, it is easy to see that the inclusion $ C^G \subset C $ induces an isomorphism
$$ C^G/I^G \cong (C/I)^G. $$
But now we compute that
$$ (C/(z))^G = k[[w]]^G = k[[w^{n''}]] = k[[t]], $$
and that
$$ (C/(w))^G = k[[z]]^G = k[[z^{n''}]] = k[[s]]. $$
So it follows that the formal branches are smooth.

In order to show that there are two irreducible components passing through $x'$, we first note that if  $x$ is the image of $x'$ in the morphism
$$ f : \mathcal{X}' \rightarrow \mathcal{X}, $$
then by assumption, there are precisely two irreducible components of the special fiber meeting at $x$. By \cite{Liubook}, Lemma 10.4.34, these lift to distinct irreducible components meeting at $x'$. But since the completion map
$$ \mathcal{O}_{\mathcal{X}',x'} \rightarrow \widehat{\mathcal{O}}_{\mathcal{X}',x'} \cong C^G, $$
is faithfully flat, it follows that there are \emph{exactly} two irreducible components meeting at $x'$.

It only remains to prove that these components are smooth at $x'$. Let $Z$ be one of the components, and denote by $ J \subset \mathcal{O}_{\mathcal{X}',x'} $ the ideal of $Z$. We have a canonical isomorphism
\begin{equation}\label{equation 1}
\widehat{\mathcal{O}}_{Z,x'} \cong \widehat{(\mathcal{O}_{\mathcal{X}',x'}/J)} \cong \widehat{\mathcal{O}}_{\mathcal{X}',x'}/ J \cdot \widehat{\mathcal{O}}_{\mathcal{X}',x'},
\end{equation}
for the completions in the various maximal ideals.

We clearly have that $ V(J \cdot \widehat{\mathcal{O}}_{\mathcal{X}',x'}) \subset \Spec(\widehat{\mathcal{O}}_{\mathcal{X}',x'}) $ is irreducible. Furthermore, since $ \mathcal{O}_{Z,x'} $ is excellent and reduced, it follows that $ \widehat{\mathcal{O}}_{Z,x'} $ is reduced (\cite{Liubook}, Proposition 8.2.41). So $V(J \cdot \widehat{\mathcal{O}}_{\mathcal{X}',x'}) \subset \Spec(\widehat{\mathcal{O}}_{\mathcal{X}',x'})$ is also reduced, and consequently $ J \cdot \widehat{\mathcal{O}}_{\mathcal{X}',x'} $ must be the ideal of one of the branches of the special fiber. Hence $ \widehat{\mathcal{O}}_{Z,x'} $ is smooth, and therefore $Z$ is smooth at $x'$.
\end{pf}

\section{The special fiber of the normalization}\label{section 3}

In the previous section we studied the local analytic structure of the morphism $ f : \mathcal{X}' \rightarrow \mathcal{X} $. In this section, we will investigate the special fiber of this map.

\subsection{}\label{section 3.1}
Let $ E \subset \mathcal{X}_k $ be an irreducible component. We make the following assumptions on $E$:

(1) $ E \cong \mathbb{P}^1_k $.

(2) $ E $ meets the rest of the special fiber at exactly two points $ \{ x_1, x_2 \} $.

(3) Let $ F_1, F_2 $ be the components meeting $E$ (where possibly $F_1 = F_2$), and let $a_i$ be the multiplicity of $F_i$. Let $b$ denote the multiplicity of $E$. Then we assume that $b$ is not divisible by $p$, or that both $a_1$ and $a_2$ are not divisible by $p$.

In other words, we demand that at least one of the components in the pair $(E,F_i)$ has multiplicity that is not divisible by $p$, for $ i = 1, 2 $. 

\subsection{} 
The following proposition is the key result in this section.

\begin{prop}\label{prop. 3.1}
Let us keep the hypotheses in Section \ref{section 3.1}. Then $ (f^{-1}(E))_{\Red} $ consists of a disjoint union of smooth and rational curves. Each of these curves meets the rest of the special fiber at exactly two points.
\end{prop}
\begin{pf}
Let $ C' := (f^{-1}(E))_{\Red} $. Then $f$ induces a finite morphism
$$ f_E: C' \rightarrow E. $$

Let $ x \in E - \{ x_1, x_2 \} $, and let $ x' $ be a point of $ \mathcal{X}'_{k} $ mapping to $x$. From the computations in Section \ref{section 2.7}, it follows that the map 
$$ \widehat{\mathcal{O}}_{E,x} \rightarrow \widehat{\mathcal{O}}_{C',x'} $$ 
of the completions of the local rings at $x$ and $x'$ can be described by 
$$ k[[u]] \rightarrow k[[s]], $$
where $u \mapsto s$. In particular we see that $ (f^{-1}(E))_{\Red} $ is regular at all points above $ E - \{ x_1, x_2 \} $, and that $ f_E : C' \rightarrow E $ is \'etale of degree $ \Deg(f_E) = \Gcd(b,n) $ above $ E - \{ x_1, x_2 \} $.

Assume now that $ x_i \in \{ x_1, x_2 \} $, and let $ d_i = \Gcd(a_i,b,n) $. From the equality 
$$ E \cdot \mathcal{X}_k = a_1 + b (E \cdot E) + a_2 = 0, $$ 
it follows easily that $ \Gcd(a_1,b) =  \Gcd(a_2,b) $ and hence $ d_1 = d_2 $. This integer will therefore be written $d$. From the computations in Section \ref{section 2.6}, it follows that 
$$ \widehat{\mathcal{O}}_{E,x} \rightarrow \widehat{\mathcal{O}}_{C',x'} $$ 
can be described by
$$ k[[u]] \rightarrow k[[s]], $$
where $u \mapsto s^c$, and where $ c = \Gcd(b,n)/d $. In particular, $ (f^{-1}(E))_{\Red} $ is regular also above $ \{ x_1, x_2 \} $, and $ f_E : C' \rightarrow E $ is tamely ramified with index $c$ at each of the $d$ points mapping to $x$.

Let $ E'_1, \ldots , E'_l $ denote the irreducible components of $C'$. Note that $ l \leq d $. The Riemann-Hurwitz formula then gives that
$$ \sum_{j=1}^l (2p_a(E'_j)-2) = \Deg(f_E)(2p_a(E)-2) + \sum_{x' \mapsto x_1} (c-1) + \sum_{x' \mapsto x_2} (c-1) $$
$$ = -2cd + d(c-1) + d(c-1) = -2d. $$
But this equation can only be fulfilled if $ l = d $ and $ p_a(E'_j) = 0 $ for all $j$.

We also see that if $E'$ is any of the components of $C'$, then there is a unique point $x'_i$ of $E'$ that maps to $x_i$, for $ i = 1, 2 $, and that $ x'_1, x'_2 $ are exactly the  points where $E'$ meets the rest of the special fiber.
\end{pf}

\subsection{}\label{section 3.3}
Let $ E \subset \mathcal{X}_k $ be an irreducible component. We make the following assumptions on $E$:

(1) $ E \cong \mathbb{P}^1_k $.

(2) $ E $ meets the rest of the special fiber at exactly one point $ \{ x_0 \} $.

(3) Let $F$ be the component meeting $E$, and let $a$ be the multiplicity of $F$. Let $b$ denote the multiplicity of $E$. Then we assume that $b$ is not divisible by $p$.

\begin{prop}\label{prop. 3.2}
Let us keep the assumptions above. Then $ (f^{-1}(E))_{\Red} $ consists of a disjoint union of smooth and rational curves. Each of these curves meets the rest of the special fiber at exactly one point.
\end{prop}
\begin{pf}
Let $ C' := (f^{-1}(E))_{\Red} $. We shall again consider the induced morphism $ f_E : C' \rightarrow E $.

Let $ x \in E - \{ x_0 \} $, and let $ x' $ be a point of $ C' $ mapping to $x$. From the computations in Section \ref{section 2.7}, we see that the map
$$ \widehat{\mathcal{O}}_{E,x} \rightarrow \widehat{\mathcal{O}}_{C',x'} $$ 
of the completions of the local rings at $x$ and $x'$ can be described by 
$$ k[[u]] \rightarrow k[[t]], $$
where $ u \mapsto t $. In particular, it follows that $ C' $ is regular above $ E - \{ x_0 \} $, and that $ f_E : C' \rightarrow E $ is \'etale of degree $ \Deg(f_E) = \Gcd(b,n) $ above $ E - \{ x_0 \} $.

From the equality $ E \cdot \mathcal{X}_k = a + b (E \cdot E) = 0 $, it follows that $ b $ divides $ a $. Hence also $ \Gcd(b,n) $ divides $a$. In particular, $ \Gcd(b,n) = \Gcd(a,b,n) $. A similar computation as in the proof of Proposition \ref{prop. 3.1} now shows that $ C' $ is regular at all points mapping to $x_0$, and that $ f_E : C' \rightarrow E $ is \'etale above $x_0$.
So the result follows.
\end{pf}

\section{Tame cyclic quotient singularities}\label{section 4}

Recall the setup from Section \ref{section 2.1}. We had an SNC-model $ \mathcal{X}/S $, and made a tamely ramified finite base change $ S' \rightarrow S $. Then we obtained a model $ \mathcal{X}' $ by normalizing the pullback of $ \mathcal{X} $ to $S'$.

Our aim in this section is to describe precisely the singular locus of $ \mathcal{X}' $ and to study the minimal desingularization $ \mathcal{X}'_{\Md} \rightarrow \mathcal{X}' $. Recall that the minimal desingularization $ \mathcal{X}'_{\Md} $ of $ \mathcal{X}' $ is the unique regular $S'$-scheme with a proper birational $S'$-morphism 
$$ \rho : \mathcal{X}'_{\Md} \rightarrow \mathcal{X}',$$ 
such that the exceptional locus of $\rho$ contains no $(-1)$-curves. The existence of a minimal desingularization follows since $ \mathcal{X}' $ is excellent. Furthermore, $ \rho : \mathcal{X}'_{\Md} \rightarrow \mathcal{X}'$ is an isomorphism outside the singular locus of $\mathcal{X}'$ (see \cite{Lip}, \cite{LipDesing}).

We will now define a class of surface singularities that will contain all singularities that $ \mathcal{X}' $ may possibly have.

\subsection{Tame cyclic quotient singularities}
Let $ \mathcal{Y} $ be a normal curve over a complete discrete valuation ring $D$ with algebraically closed residue field $k$. Let us assume that $\mathcal{Y}$ has smooth generic fiber, so that $\mathcal{Y}_{\Sing}$ consists of finitely many closed points in the special fiber $\mathcal{Y}_k $.

\begin{dfn}
A closed point $y$ in the special fiber $\mathcal{Y}_{k}$ is a \emph{tame cyclic quotient singularity} if there exists a positive integer $n > 1$ not divisible by $ p = \Char (k) $, a unit $ r \in (\mathbb{Z}/n \mathbb{Z})^* $, and integers $ m_1 > 0 $ and $ m_2 \geq 0 $ satisfying $ m_1 \equiv - r m_2 $ mod $n$ such that $ \widehat{\mathcal{O}_{\mathcal{Y},y}} $ is isomorphic to the subalgebra of $ \boldsymbol{\mu}_{n}(k) $-invariants in 
$$ D[[t_1,t_2]]/(t_1^{m_1} t_2^{m_2} - \pi), $$
under the action given by $ [\xi](t_1) = \xi t_1 $ and $ [\xi](t_2) = \xi^r t_2 $, for any $ \xi \in \boldsymbol{\mu}_{n}(k) $, where $ \pi $ is a uniformizing parameter for $D$.
\end{dfn}

\begin{rmk}
This is a simplified version of Definition 2.3.6 in \cite{CED}. Even though the definition involves a choice of coordinates, one can actually show that the parameters $n$ and $r$ are \emph{intrinsic} for the singularity $y$ (\cite{CED}, Remark 2.3.7). The case of one analytic branch, $ m_2 = 0 $, will not occur in this paper.
\end{rmk}

\subsection{The singular locus of $ \mathcal{X}' $}
Let us now describe the singular locus of $ \mathcal{X}' $. Note that since $ \mathcal{X}' $ is normal, and the generic fiber is smooth, $ \mathcal{X}'_{\Sing} $ consists of finitely many closed points in the special fiber $ \mathcal{X}'_{k} $.

Let $ x' \in \mathcal{X}' $ be a closed point in the special fiber, and let $ x = f(x') $ be the image of $x'$ under $ f : \mathcal{X}'  \rightarrow \mathcal{X} $. Then $x$ is a closed point in $ \mathcal{X}_s $, and the local analytic structure of $ \mathcal{X}' $ at $x'$ depends only on the local analytic structure of $ \mathcal{X} $ at $x$.

\begin{prop} Let $\mathcal{X}'$ be as above, and let $ x' \in \mathcal{X}'_k $ be a closed point. 
\begin{enumerate}
\item Assume $ x = f(x') $ belongs to a unique irreducible component of $ \mathcal{X}_k $. Then $ \mathcal{X}' $ is regular at $x'$.

\item Assume that $ x = f(x') $ is an intersection point of two distinct irreducible components with multiplicities $a$ and $b$. Let $ a'' = a/\Gcd(a,n) $, $ b'' = b/\Gcd(b,n) $ and $ n'' = n \cdot \Gcd(a,b,n)/(\Gcd(a,n) \cdot \Gcd(b,n)) $. Let $ r \in (\mathbb{Z}/n'')^* $ be such that $ rb'' + a'' \equiv 0 $ $ \Mod $ $ n''$. Then $ \mathcal{X}' $ has a tame cyclic quotient singularity at $x'$, formally isomorphic to the invariant ring 
$$ \{R'[[z,w]]/(\pi' - z^{a''} w^{b''})\}^{\boldsymbol{\mu}_{n''}(k)}, $$
where the $ \boldsymbol{\mu}_{n''}(k) $-action on $ R'[[z,w]]/(\pi' - z^{a''} w^{b''}) $ is given by $ [\xi](z) = \xi z $ and $ [\xi](w) = \xi^r w $, for any $ \xi \in \boldsymbol{\mu}_{n''}(k) $. 
\end{enumerate}
\end{prop}
\begin{pf}
Part (1) follows from Proposition \ref{prop. 2.8}, since $ \widehat{\mathcal{O}}_{\mathcal{X}',x'} $ is regular if and only if $ \mathcal{O}_{\mathcal{X}',x'} $ is regular (\cite{Liubook}, Lemma 4.2.26). Part (2) follows from Proposition \ref{prop. 2.5}. 
\end{pf}

\subsection{Minimal resolution of tame cyclic quotient singularities}
Let $ y \in \mathcal{Y} $ be a tame cyclic quotient singularity, with parameters $m_1$, $m_2$, $n$ and $r$. The minimal desingularization $ \rho : \mathcal{Y}_{\Md} \rightarrow \mathcal{Y} $ locally around $y$ is completely described in Theorem 2.4.1 in \cite{CED}. The properties that we shall need are listed below.
\begin{enumerate}
\item The fiber of $\rho$ over $y$ consists of a chain $ E_1, \ldots ,E_{\lambda} $.

\item All $E_j$ are isomorphic to $ \mathbb{P}^1_{k} $.

\item All intersections in this chain are transverse.

\item $E_1$ is transverse to the strict transform $ \widetilde{\mathcal{Y}}_1 $ of the component $\mathcal{Y}_1$ of $\mathcal{Y}_k$ through $y$ with multiplicity $m_2$, and similarly for $E_{\lambda}$ and the component $ \widetilde{\mathcal{Y}}_2 $ with multiplicity $m_1$.

\item The self intersection numbers of the $E_j$ are determined in terms of the parameters of the singularity. Consider the Jung-Hirzebruch continued fraction expansion
$$ \frac{n}{r} = b_1 -  \frac{1}{b_2 - \frac{1}{\ldots - \frac{1}{ b_{\lambda}}}}. $$
Then we have that $ E_j \cdot E_j = - b_j < -1 $ for all $j$.
\end{enumerate}

The multiplicity $ \mu_j $ in $ (\mathcal{Y}_{\Md})_k $ of the component $ E_j $ is determined in the following fashion: The special fiber of $ \mathcal{Y}_{\Md} $ (as a divisor) is 
$$ (\mathcal{Y}_{\Md})_k = m_2 \widetilde{\mathcal{Y}}_1 + \mu_1 E_1 + \mu_2 E_2 + \ldots + \mu_{\lambda} E_{\lambda} + m_1 \widetilde{\mathcal{Y}}_2 + \ldots $$
Since $ \mathcal{Y}_{\Md} $ is an arithmetic surface, we have that $ E_j \cdot (\mathcal{Y}_{\Md})_k = 0$ (\cite{Liubook}, Proposition 9.1.21). This gives the equation $ \mu_{j-1} + \mu_{j+1} - b_j \mu_j = 0 $, where we put $ \mu_0 = m_2 $ and $ \mu_{\lambda + 1} = m_1 $. Doing this for every $j$ gives a system of equations that can be solved for the $\mu_j$-s.

\begin{cor}
The minimal desingularization $ \mathcal{X}'_{\Md} $ of $ \mathcal{X}' $ is an SNC-model.
\end{cor}
\begin{pf}
We saw in Corollary \ref{cor. 2.9} and Proposition \ref{prop. 2.6} that the irreducible components of $ \mathcal{X}'_k $ are smooth. In particular, $ \mathcal{X}'_k $ is a strict normal crossings divisor when we restrict to the regular locus of $ \mathcal{X}' $. Furthermore, $ \mathcal{X}' $ has tame cyclic quotient singularities only at closed points in the special fiber where two distinct irreducible components meet. Since the minimal desingularization is an isomorphism outside the singular locus, it follows from the description above that $ \mathcal{X}'_{\Md} $ is indeed an SNC-model.
\end{pf}

\section{Contraction to semi-stability}

In this section we prove a technical lemma which will be important in Sections 6 and 7.

\begin{lemma}\label{lemma 5.1}
Let $ \mathcal{Y}/S $ be a regular surface with normal crossings, where $S$ is the spectrum of a discrete valuation ring, with algebraically closed residue field $k$. Assume that a part of the special fiber $ \mathcal{Y}_k $ consists of irreducible components 
$$ \mathcal{E} : E_0, E_1, \ldots, E_l, E_{l+1}, $$ 
such that
\begin{enumerate}
\item $E_i \cong \mathbb{P}^1_k $ for all $ i \in \{1, \ldots ,l \} $.

\item $ E_i \cap E_{i+1} $ is a unique point for all $ i \in \{0, \ldots ,l \} $, except for the case where $l=1$ and $ E_0 = E_{l+1} $, where we assume that $E_1$ meets $E_0$ in exactly two distinct points.

\item $E_i$ meets no components of $ \mathcal{Y}_k $ other than $E_{i-1}$ and $E_{i+1}$, for all $ i \in \{1, \ldots ,l \} $.

\item $\Mult(E_0) = \Mult(E_{l+1}) = 1$.
\end{enumerate} 

Then there exists a proper birational morphism $ \Psi : \mathcal{Y} \rightarrow \mathcal{Z} $, relative to $S$, such that $ \mathcal{Z} $ is a regular surface with normal crossings, and if $ E_i \in \mathcal{E} $ is a component such that $ \Mult(E_i) > 1 $, then $ E_i $ is contracted under $ \Psi $.
\end{lemma}
\begin{pf}
We shall construct $ \Psi $ stepwise, by contracting one component at the time. Let $ m_i = \Mult(E_i) $, for $ i \in \{1, \ldots, l \} $. Assume that $ m_i \in \{ m_1 , \ldots, m_l \} \subset \mathbb{N} $ is maximal, and that $ m_i > 1 $ (otherwise we are done). Then we necessarily have that $ m_i > m_{i+1} $ and $ m_i > m_{i-1} $. For if it was not so, then $ m_i = m_{i+1} $ or $ m_i = m_{i-1} $. Let us assume $ m_i = m_{i+1} $. Intersecting the special fiber $ \mathcal{Y}_k $ with $E_i$ would then give the equation 
$$ 0 = E_i \cdot \mathcal{Y}_k = E_i \cdot ( \sum_{j=0}^{l+1} m_j E_j + \ldots ) = m_{i-1} + m_i E_i^2 + m_{i+1}. $$
From this equation it follows that $ m_i $ also equals $ m_{i-1} $. Continuing inductively, it is easy to see that $ m_i $ equals $ m_j $ for \emph{all} $ j \in \{0, \ldots ,l+1 \} $. But by assumption $ m_0 = m_{l+1} = 1 $, so this is a contradiction. Hence, if $ m_i $ is maximal in the set of multiplicities, and $ m_i > 1 $, then $ m_i > m_{i+1} $ and $ m_i > m_{i-1} $.

Next, we show that a component $E_i$ with maximal multiplicity is necessarily a $(-1)$-curve. Indeed, from the equality
$$ m_{i-1} + m_i E_i^2 + m_{i+1} = 0, $$ 
it follows that
$$ E_i^2 = -(m_{i-1} + m_{i+1})/m_i > - (m_i + m_i)/m_i = - 2. $$
But from the above equation it is clear that $ E_i^2 < 0 $, so $ E_i^2 = -1 $. Since $ E_i $ is smooth and rational, it may therefore be contracted.

Let $ \Psi_1 : \mathcal{Y} \rightarrow \mathcal{Z}_1 $ be the contraction morphism of $E_i$. We claim that $ \mathcal{Z}_1 $ is a regular fibered surface with normal crossings. To see this, let us first assume that $ l > 1 $, or that $E_0$ and $ E_{l+1} $ are distinct components. Then the statement follows from \cite{Liubook}, Lemma 9.3.35. In the case where $ l = 1 $ and $ C:= E_0 = E_{l+1} $, we have that $ \Mult(E_1) = 2 $, and hence the image of $C$ has multiplicity $2$ in the image point $Q$ of $E_1$ on $ \mathcal{Z}_1 $. It follows that $C$ has a nodal singularity at $Q$, and hence $ \mathcal{Z}_1 $ has normal crossings at $Q$.

The image of $ \mathcal{E} $ on $ \mathcal{Z}_1 $ is a new chain 
$$ \mathcal{E}_1 : E_0^1, E_1^1, \ldots, E_{l-1}^1, E_{l}^1, $$
that satisfies the conditions (1)-(4) above. So we may choose a maximal element
$$ \mu_{i}^1 \in \{ \mu_1^1, \ldots , \mu_{l-1}^1 \} \subset \mathbb{N}, $$ 
where $ \mu_j^1 $ denotes the multiplicity of $E_j^1 $. If $ \mu_{i}^1 > 1 $, the component $ E_{i}^1 $ will have self intersection $-1$ by the same argument as above, and we therefore get a contraction morphism
$$ \Psi_2 : \mathcal{Z}_1 \rightarrow \mathcal{Z}_2, $$
where $ \mathcal{Z}_2 $ is a regular fibered surface with normal crossings, and $ E_{i}^1 $ is contracted to a point. Continuing inductively, we get a series of contractions
$$ \mathcal{Y} \rightarrow \mathcal{Z}_1 \rightarrow \ldots \rightarrow \mathcal{Z}_r, $$
where the image in $ \mathcal{Z}_r $ of the original chain $ \mathcal{E} $ is semistable. Therefore we let $ \mathcal{Z} = \mathcal{Z}_r $, and take $ \Psi $ to be the composition of all the $ \Psi_i $.
\end{pf}

\section{chains, principal components}

\subsection{Chains} 

To begin with, we only assume that $ \mathcal{X}/S $ is a normal fibered surface. 

\begin{dfn}
A \emph{chain} in $ \mathcal{X}_k $ is a set of smooth and rational curves $ E_1, E_2, \ldots, E_l $ such that all the $E_i$ intersect the rest of $ \mathcal{X}_k $ in at most two points, and $E_i$ meets only one component at each point where it intersects with the rest of $ \mathcal{X}_k $. Furthermore, $E_i$ intersects $E_{i+1}$ in a unique point for $ i = 1, \ldots, l-1 $. We call $E_1$ and $E_l$ the \emph{ends} of the chain. 

If $E$ is another smooth and rational component of $ \mathcal{X}_k $ meeting the rest of the special fiber in at most two points, and meeting one of the ends of the chain, we may extend the chain by $E$. If the chain cannot be extended in this way, we say that it is \emph{maximal}. A chain where the two ends meet the other components of the chain twice is called a \emph{loop}. Otherwise, we say that the chain is \emph{open}.
\end{dfn}

\subsection{Principal components}

Following \cite{Xiao}, we make the following definition.
\begin{dfn}\label{principal}
Let $ \mathcal{X}/S $ be a regular model with normal crossings. We say that an irreducible component $F$ of $ \mathcal{X}_k $ is \emph{principal} if either
\begin{enumerate}
\item $F$ is smooth and rational, and meets the rest of the special fiber at \emph{more} than two points, or

\item $F$ is \emph{not} smooth and rational
\end{enumerate}
\end{dfn}

Hence, if $E$ is a non-principal irreducible component, then $E$ is smooth and rational, and meets the rest of the special fiber in at most two points.

Let $ \mathcal{E} =\{ E_1, E_2, \ldots, E_l \} $ be a chain in $ \mathcal{X}_k$. If $E_1$ (resp. $E_l$) meets an ``outer'' component, this will be denoted by $E_0$ (resp. $E_{l+1}$). Let $a_i$ be the multiplicity of $E_i$, for $ i = 0, 1, \ldots, l, l+1 $. We will allways assume that $p$ does not divide both $a_i$ and $a_{i+1}$ for all $i$.

Let $ S' \rightarrow S $ be a tamely ramified extension of degree $n$, and $ f : \mathcal{X}' \rightarrow \mathcal{X} $ the induced morphism of models as in Section \ref{section 2.1}. We shall now describe the inverse image of the chain $ \mathcal{E} $ under $f$. Notice that this is a generalization of the situation studied in Section \ref{section 3}.

\begin{lemma}\label{lemma 6.3}
Let us keep the hypotheses above. Let the integers $d_i$ be defined by $ d_i = \Gcd(a_i,a_{i+1},n) $, where $i = 0, 1, \ldots , l $. Then we have that the $d_i$ are equal for all $i$. We denote this integer by $d_{\mathcal{E}}$. 
\end{lemma}
\begin{pf}
By intersecting $E_i$ with the special fiber, we get the equation
$$ 0 = E_i \cdot \mathcal{X}_k = E_i \cdot (a_{i-1} E_{i-1} + a_i E_i + a_{i+1} E_{i+1} + \ldots ) = a_{i-1} + a_i E_i^2 + a_{i+1}, $$
from which it easily follows that $ \Gcd(a_{i-1},a_i) = \Gcd(a_i,a_{i+1}) $ for all $i = 1, \ldots , l $. By taking the greatest common divisor with $n$, we get that all the integers $d_i$ indeed are equal.
\end{pf}

\begin{prop}\label{lemma 6.4}
Let us keep the hypotheses above. Let $ \mathcal{E} =\{ E_1, \ldots, E_l \} $ be a chain that is not a loop, and assume that at least one of the components $E_1$ and $ E_{l} $ meets an outer component. Then the reduced inverse image of the chain $ \mathcal{E} =\{ E_1, \ldots, E_l \} $ is a disjoint union of $d_{\mathcal{E}}$ chains of smooth and rational curves.
\end{prop}
\begin{pf}
From Proposition \ref{prop. 3.1}, Proposition \ref{prop. 3.2} and Lemma \ref{lemma 6.3}, we know that the inverse image of each $E_i$ is a disjoint union of $d_{\mathcal{E}}$ smooth and rational curves $E'_{ij}$. Each $E'_{ij}$ meets the rest of the special fiber $ \mathcal{X}'_k $ in exactly the same number of points as $E_i$ meets the rest of $ \mathcal{X}_k $. Furthermore, each component of $ \mathcal{X}_k $ that meets $E_i$ can be lifted to a component that meets $E'_{ij}$. The result follows immediately.
\end{pf}

\begin{rmk}\label{remark 6.5}
Let $ \mathcal{X}'_{\Md} \rightarrow \mathcal{X}' $ be the minimal desingularization. If $ \mathcal{E}' $ is a chain (or loop) of smooth and rational curves in the special fiber of $ \mathcal{X}' $, then it follows from the description of the minimal desingularization in Section \ref{section 4} that the inverse image $ \bar{\mathcal{E}}' $ of $ \mathcal{E}' $ on $ \mathcal{X}'_{\Md} $ is again a chain (or loop) of smooth and rational curves. Indeed, $ \mathcal{X}' $ has only tame cyclic quotient singularities, and these are located at closed points where two branches of the special fiber meet. The inverse image of a singularity consists of a chain of smooth rational curves. So $ \bar{\mathcal{E}}' $ is just a ``blowup'' of $ \mathcal{E}' $.
\end{rmk}

\subsection{}
We shall now consider the case where $ \mathcal{X}/S $ is a minimal SNC-model satisfying Saito's criterion, and where the generic fiber has genus at least equal to $2$.

Let $ \mathcal{E} = \{ E_1, E_2, \ldots, E_l \}$ be a maximal chain of components in the special fiber ~$ \mathcal{X}_k $. The following result shows that the chain must meet a principal component.

\begin{prop}\label{lemma 6.6}
Let $ \mathcal{E} = \{ E_1, E_2, \ldots, E_l \}$ be a maximal chain of smooth and rational curves in $ \mathcal{X}_k $. Then $ \mathcal{E} $ cannot be a loop. Furthermore, at least one of the ends $E_1$ or $E_l$ of the chain must meet a principal component.
\end{prop}
\begin{pf}
Let us assume that $ \mathcal{E} $ does not meet a principal component. In that case it is easy to see that $ \mathcal{E} $ equals the whole special fiber (since the special fiber is connected), and that $ \mathcal{E} $ is either (a) a loop, or (b) an open chain.

We shall first show that case (a) leads to a contradiction. We note that at least one of the components must have multiplicity not divisible by $p$. Indeed, since $ \mathcal{X} $ has strict normal crossings, all components are smooth, so the loop must consist of at least two components. But then all components cannot have multiplicity divisible by $p$, since that would contradict the assumption that $ \mathcal{X} $ fulfills Saito's criterion.

Let $a_i$ be the multiplicity of $E_i$, and let $ \bar{a}_i $ be the prime-to-${p}$ part of $a_i$. Let $ n = \Lcm_{1 \leq i \leq l} \{ \bar{a}_i \} $, and let $ S' \rightarrow S $ be a tamely ramified extension of degree $n$. Arguing as in Proposition \ref{lemma 6.4}, it follows that the inverse image of $ \mathcal{X}_k $ is a disjoint union of loops. Since $ \mathcal{X}'_{k} $ is connected, we get that $ \mathcal{X}'_{k} $ is a loop. Let $ \mathcal{X}'_{\Md} \rightarrow \mathcal{X}' $ be the minimal desingularization. Then it follows from Remark \ref{remark 6.5} that the special fiber of $ \mathcal{X}'_{\Md} $ is a loop. But we also have that some of the irreducible components are reduced, corresponding to the $a_i$ that are not divisible by $p$. By Lemma \ref{lemma 5.1}, we may then contract all components with multiplicity greater than one. But this means that the minimal regular model $ \mathcal{X}'_{\Min} $ is semi-stable, and that the special fiber is either a rational curve with exactly one node, or a loop of smooth and rational curves. In any case, we see that $ p_a((\mathcal{X}'_{\Min})_k) = 1 $, which is a contradiction.

In case (b), a similar argument implies $ p_a((\mathcal{X}'_{\Min})_k) = 0 $, again a contradiction.
\end{pf}

\section{The minimal extension realizing stable reduction}

Let $R$ be a complete discrete valuation ring with algebraically closed residue field $k$, and let $X/K$ be a smooth, projective, geometrically connected curve of genus at least equal to two. Let $ \mathcal{X}/S $ be the minimal SNC-model of $X$, and assume that $ \mathcal{X}/S $ satisfies Saito's criterion.

Denote by $ \mathcal{F} $ the set of irreducible components of $ \mathcal{X}_k $ that are principal, and let $ n = \Lcm \{ \Mult(F) | F \in \mathcal{F} \} $.

\begin{thm}\label{theorem 7.1}
Let $S' \rightarrow S$ be a tamely ramified extension of degree $n$. Then $ X_{K'} $ has stable reduction over $S'$.
\end{thm}
\begin{pf}
Let $ f : \mathcal{X}' \rightarrow \mathcal{X} $ be the morphism constructed in Section \ref{section 2.1}, and let $ \rho : \mathcal{X}'_{\Md} \rightarrow \mathcal{X}' $ be the minimal desingularization. Consider the open subscheme $ \mathcal{V} = \mathcal{X} - \cup_{ E \notin \mathcal{F} } E $ of $ \mathcal{X} $, and let $ \mathcal{V}' = f^{-1}(\mathcal{V}) $, and furthermore $ \mathcal{U}' = \rho^{-1}(\mathcal{V}') $. Let us first show that $ \mathcal{U}' $ is semi-stable over $S'$. Indeed, if $ F \in \mathcal{F} $, and $ m = \Mult(F) $, then we have that $ m | n $. Consequently, if $ x' \in \mathcal{X}'_{k} $ is a closed point, and $ x' \in \mathcal{V}' $, then it follows from the local computations in Section \ref{section 2} that we either have
$$ \widehat{\mathcal{O}}_{\mathcal{X}',x'} \cong R'[[u,v]]/(\pi'^{n'}-v), $$
or that we have
$$ \widehat{\mathcal{O}}_{\mathcal{X}',x'} \cong R'[[u,v]]/(\pi'^{n''}-uv). $$
So the inverse image of $x'$ under $ \rho : \mathcal{X}'_{\Md} \rightarrow \mathcal{X}' $ can at most be a chain of \emph{reduced} smooth and rational curves, where all intersections are transversal (cf. ~\cite{DelMum}, Proposition 2.3). It follows that $ (\mathcal{U}')_{k} $ is reduced and has only nodal singularities.

Now, let $ \mathcal{E} =  \{ E_1, E_2, \ldots, E_l \} $ be a maximal chain of smooth and rational curves in $ \mathcal{X}_k $. By Proposition \ref{lemma 6.6}, we may assume that $E_1$ meets a principal component $F$. If $E_l$ does not meet a principal component, then the multiplicity of $E_l$ divides the multiplicity of $F$. Let $ \bar{\mathcal{E}}' $ be the inverse image of $ \mathcal{E} $ on $ \mathcal{X}'_{\Md} $. By Lemma \ref{lemma 6.4} and Remark \ref{remark 6.5}, we have that $ \bar{\mathcal{E}}' $ is a disjoint union of chains of smooth and rational curves. Furthermore, we always have that one of the ends of each chain meets a reduced component of the special fiber of $ \mathcal{X}'_{\Md} $, and that the other end either meets a reduced component, or is itself reduced.

We now have that $ \mathcal{X}'_{\Md} $ has strict normal crossings, and that the only components that are not reduced are parts of chains of smooth and rational curves where the ends meet reduced components. We can now apply Lemma \ref{lemma 5.1}, and get that there exists a contraction
$$ \mathcal{X}'_{\Md} \rightarrow \mathcal{X}'_0, $$
such that $ \mathcal{X}'_0 $ is regular and a semi-stable model of $X_{K'}$. It then follows from \cite{Liubook}, Theorem 10.3.34, that $X_{K'}$ has stable reduction over $S'$. 
\end{pf}

\subsection{}
Having established Theorem \ref{theorem 7.1}, we can weaken the hypothesis on the base slightly, allowing the discrete valuation ring to be only strictly henselian, which is perhaps more natural in the light of Saito's result.

\begin{cor}\label{corollary 7.2}
Let $R$ be a strictly henselian discrete valuation ring, with algebraically closed residue field $k$. Let $X$ be a smooth, projective and geometrically connected curve of genus $ g(X) \geq 2 $ over $K$, where $ K $ is the quotient field of $R$. Let $ \mathcal{X}/S $ be the minimal SNC-model of $X$, and assume that it satisfies Saito's criterion.

Let $n$ be the least common multiple of the multiplicities of the principal components of $\mathcal{X}_k$, and let $K' = K[\pi']/(\pi'^n - \pi)$. Then $ X \otimes_K K' $ has stable reduction over $R'$, where $R'$ is the integral closure of $R$ in $K'$.
\end{cor}
\begin{pf}
Let us first note that $R'$ is again a strictly henselian discrete valuation ring, with uniformizing parameter $\pi'$ (cf. ~\cite{Neukirch}, Proposition II.7.7 and Proposition II.8.2). Since $R'$ is finite over $R$, we have that
$$ \widehat{R'} = R' \otimes_R \widehat{R}, $$
that $ \pi $ (resp. $ \pi' $) is a uniformizing parameter of $\widehat{R} $ (resp. $\widehat{R'}$), and that $$ \widehat{R'} = \widehat{R}[\pi']/(\pi'^n - \pi). $$
The pullback 
$$ \widehat{\mathcal{X}} = \mathcal{X} \times_{\Spec(R)} \Spec(\widehat{R}) $$
is the minimal SNC-model of $ \widehat{X} = X \otimes_K \widehat{K} $. To see this, let us first show that $ \widehat{\mathcal{X}} $ is regular. The generic fiber is $ \widehat{X} $, which is smooth, so we need only check points in the special fiber. But the projection morphism
$$ \widehat{\mathcal{X}} \rightarrow \mathcal{X} $$
induces an isomorphism of the special fibers, and then it follows by \cite{Liubook}, Lemma 8.3.49 (b), that $ \widehat{\mathcal{X}} $ is regular. Furthermore, as the special fibers are isomorphic, we get indeed that $ \widehat{\mathcal{X}} $ is the minimal SNC-model. 

By Theorem \ref{theorem 7.1}, the extension $ \widehat{R} \rightarrow \widehat{R'} $ will realize stable reduction for $ \widehat{X} $. Let now $ \widehat{\mathcal{Y}} $ be the canonical model of $ \widehat{X'} $ over $ \Spec(\widehat{R'}) $, and let $ \mathcal{Y} $ be the canonical model of $X' $ over $ \Spec(R') $. By \cite{Liubook}, Proposition 10.1.17, we have a canonical isomorphism
$$ \widehat{\mathcal{Y}} \cong \mathcal{Y} \times_{\Spec(R')} \Spec(\widehat{R'}), $$
inducing an isomorphism on the special fibers of $ \widehat{\mathcal{Y}} $ and $ \mathcal{Y} $, so consequently we get that $ \mathcal{Y} $ is stable over $ \Spec(R') $.  
\end{pf}

\subsection{Minimality}
We will now prove that the extension found in Theorem \ref{theorem 7.1} is the minimal extension that realizes stable reduction for $X/K$. In order to do this, we will use a certain quotient construction, due to E. Viehweg (\cite{Viehweg}), and generalized by D. Lorenzini (\cite{Dino}). We refer to these papers for more details.

\subsection{}
Let $ K \subset L $ be the minimal extension realizing stable reduction for $X/K$. This extension is tamely ramified of some degree $ d = [L:K] $, and Galois with cyclic group $G$.

Let $ \sigma \in G $ be a generator of the Galois group, and denote also by $ \sigma $ the induced automorphism $ \sigma : \Spec(L) \rightarrow \Spec(L) $. We have that $\sigma$ induces an automorphism of $ X_{L} = X \times_{\Spec(K)} \Spec(L) $ by
$$ \Id \times \sigma : X \times_{\Spec(K)} \Spec(L) \rightarrow X \times_{\Spec(K)} \Spec(L), $$
and hence $ G $ acts on $ X_{L} $.

Let $ \mathcal{Y}  $ be the minimal regular model of $ X_{L} $ over $S_L := \Spec(R_L)$, where $R_L$ denotes the integral closure of $R$ in $L$. Since $ X_{L} $ has stable reduction over $S_L$, we have that $ \mathcal{Y}/S_L $ is a semi-stable model. By the universal property of the minimal regular model, any automorphism of the generic fiber extends uniquely to $ \mathcal{Y} $, so there exists a unique automorphism $\tau$ of $ \mathcal{Y} $ making the following diagram commute:
$$ \xymatrix{
 \mathcal{Y} \ar[d] \ar[r]^{\tau} & \mathcal{Y} \ar[d]\\
S_L \ar[r]^{\sigma} & S_L. } $$
Hence the $G$-action on $ X_{L} $ extends to $ \mathcal{Y} $.

One can now form the quotient of $ \mathcal{Y} $ for the action of $G$. Let us denote this quotient by $ \mathcal{Z} $, and let $ g : \mathcal{Y} \rightarrow \mathcal{Z} $ be the quotient morphism. Below follow some useful properties of $ \mathcal{Z} $ and $g$.

\begin{prop}\label{proposition 7.3}
(i) The quotient $ \mathcal{Z} $ is a fibered surface over $S$, with generic fiber equal to $X/K$. 

(ii) For any irreducible component $D$ of $ \mathcal{Y}_{k} $, let $ I_D = \{ \mu \in G ~|~ \mu_{|D} = \Id \} $, and let $Z = g(D)$. Then the multiplicity of $Z$ in $ \mathcal{Z}_k $ equals $ d/|I_D| $. \end{prop}
\begin{pf}
Let us explain (ii). Let $ \tau_{k} $ denote the induced automorphism on the special fiber $ \mathcal{Y}_{k} $. Since the extension $ K \subset L $ was minimal with the property of realizing stable reduction for $X/K$, we get, by Lemma 3.4 in \cite{Viehweg}, that $ \Ord(\tau_{k}) = \Ord(\sigma) = d $. Since $ \mathcal{Y}/S_L $ is semi-stable, it follows from Fact IV in \cite{Dino} that $\Mult(Z) = d/|I_D| $.
\end{pf}

\vspace{0.5cm}

It is explained in \cite{Dino} and \cite{Viehweg} how one may find an explicit desingularization of $ \mathcal{Z}/S $. 

\begin{prop}\label{proposition 7.4}
There exists a regular model $ \widetilde{\mathcal{Z}}/S $ with normal crossings, and a proper birational morphism
$$ \pi : \widetilde{\mathcal{Z}} \rightarrow \mathcal{Z}, $$
such that $\pi$ is an isomorphism over the regular locus of $ \mathcal{Z} $.

The exceptional locus of $\pi$ can be described as follows; if $ z \in \mathcal{Z}_{\Sing} $, then $ \pi^{-1}(z) $ consists of a chain of smooth rational curves meeting the rest of the special fiber at at most two points corresponding to  the ends of the chain.  
\end{prop}
\begin{pf}
This is Fact V in \cite{Dino}.
\end{pf}

\vspace{0.5cm}

We can now formulate the second main result in this paper, showing that the extension in Theorem \ref{theorem 7.1} is minimal:

\begin{thm}\label{theorem 7.5}
The extension in Theorem \ref{theorem 7.1} is the minimal extension that realizes stable reduction for $X/K$.
\end{thm}
\begin{pf}
Let $F \in \mathcal{X}_k $ be a principal component, and $ m = \Mult(F) $. We need to show that $m$ divides $d$, the degree of $ K \subset L $, where $L$ is the minimal extension realizing stable reduction.

Since $ \widetilde{\mathcal{Z}} $ is regular with normal crossings, the irreducible components have at most nodal singularities. By blowing up $ \widetilde{\mathcal{Z}} $ in these points, we obtain a surjective birational morphism $ \overline{\mathcal{Z}} \rightarrow \widetilde{\mathcal{Z}} $ such that $ \overline{\mathcal{Z}} $ is an SNC-model, and such that the exceptional locus consists of smooth and rational curves meeting the rest of the special fiber in exactly two points. Furthermore, $ \overline{\mathcal{Z}} $ dominates $ \mathcal{X} $. That is, we have a birational and surjective morphism $ \rho: \overline{\mathcal{Z}} \rightarrow \mathcal{X} $, since $ \mathcal{X} $ was the minimal SNC-model of $X$. 

Let $\overline{F}$ be the strict transform of $F$ under this map. Since $ \rho $ can be factored as a series of blow-ups of closed points in the special fibers, we have that $\overline{F}$ is a principal component of $ \overline{\mathcal{Z}}_k $. So we have either that $p _a(\overline{F}) > 0 $, or $\overline{F} \cong \mathbb{P}^1_k $, and meets the rest of the special fiber in at least three points. In any case, we see that $\overline{F}$ does not belong to the exceptional locus of the composition of $ \pi : \widetilde{\mathcal{Z}} \rightarrow \mathcal{Z} $ with $ \overline{\mathcal{Z}} \rightarrow \widetilde{\mathcal{Z}} $. Consider the image $ Z $ of $\overline{F}$ in $ \mathcal{Z} $, and let $ D $ be a component of $ \mathcal{Y}_k $ mapping to $Z$. Then Proposition \ref{proposition 7.3} gives that $ d = m \cdot |I_D| $, since $ m = \Mult(F) = \Mult(Z) $. Hence $ \Mult(F) $ divides $d$, which is what we wanted to show.
\end{pf}

\subsection{}
As a Corollary, we obtain a new, more elementary proof of Theorem \ref{S-crit}. 

\begin{pf}
The sufficiency of $(*)$ in Theorem \ref{S-crit} follows from Theorem \ref{theorem 7.1}.

It remains to show the necessity of $(*)$. This is stated in \cite{Dino}, Corollary 1.1 without proof, so we give the argument here for completeness. 

Let us assume that the minimal extension $ K \subset L $ realizing stable reduction for $X/K$ is tamely ramified. Let $ \mathcal{Y}/S_L $ be the minimal regular model of $X \otimes_K L $, and consider the quotient map 
$$ g : \mathcal{Y} \rightarrow \mathcal{Y}/G = \mathcal{Z}, $$
as constructed above. Let $Z$ be any irreducible component of the special fiber of $ \mathcal{Z} $. Then we have that the multiplicity of $Z$ divides $ d = [L:K] $, and hence is prime to $p$. Let $ \widetilde{\mathcal{Z}} \rightarrow \mathcal{Z} $ be the desingularization as in Proposition \ref{proposition 7.4}. Then any irreducible component $E$ of the special fiber of $ \widetilde{\mathcal{Z}} $ with multiplicity divisible by $p$ must be exceptional. 

Assume in this case that $E$ meets the rest of the special fiber at only one point. This means that $E$ is the end of a chain, and we have seen that the multiplicity of $E$, and hence $p$, will divide the multiplicity of all components in the chain, as well as the component meeting the other end of the chain. But the other end of the chain will meet the strict transform of a component of $ \mathcal{Z}_k $, which has multiplicity prime to $p$, a contradiction.

Assume that $E$ meets the rest of the special fiber at two distinct points. Then the assumption that $E$ meets a component with multiplicity divisible by $p$ leads to a similar contradiction as above, since the chain of exceptional curves containing $E$ must eventually meet the strict transform of some irreducible component of $ \mathcal{Z}_k $. And if two successive curves in the chain have multiplicities divisible by $p$, then $p$ also divides the multiplicities of the components meeting the ends of the chain.

Consider now the morphism $ \overline{\mathcal{Z}} \rightarrow \widetilde{\mathcal{Z}} $ constructed in the proof of Theorem \ref{theorem 7.5}. Any exceptional component for this map arises as the exceptional curve of a single blow up of a closed point corresponding to a nodal singularity of the strict transform of a component of $ \mathcal{Z}_k $. Hence it is immediate that if such a curve has multiplicity divisible by $p$, it will satisfy the condition $(*)$ in Theorem \ref{S-crit}. We have now established that $ \overline{\mathcal{Z}} $ is an SNC-model, that satisfies Saito's criterion.

Finally, we need to show that $ \mathcal{X} $ fulfills Saito's criterion. So we consider $ \rho: \overline{\mathcal{Z}} \rightarrow \mathcal{X} $. This morphism can be factored as a series of blow-ups of closed points in the special fiber. Assume that some principal component of $ \mathcal{X}_k $ has multiplicity divisible by $p$. Since the strict transform of this component in $ \overline{\mathcal{Z}} $ is principal, we get a contradiction. Let $F$ be a smooth and rational component of $ \mathcal{X}_k $, with multiplicity divisible by $p$, meeting the rest of the special fiber in at most two points. Assume that this component violates Saito's criterion for $ \mathcal{X}_k $. Then the strict transform of $F$ in $ \overline{\mathcal{Z}} $ violates Saito's criterion for $ \overline{\mathcal{Z}}_k $. This follows from considering the multiplicity of the exceptional divisor $E$ when we blow up $ \mathcal{X} $ in a closed point on $F$. The multiplicity of $E$ is divisible by $p$ if we blow up in a point that does not lie on any other component, or if it is an intersection point of $F$ and another component $F'$, whose multiplicity is also divisible by $p$.  Using this argument in the series of blow-ups, it is easily seen that the strict transform of $F$ in $ \overline{\mathcal{Z}} $ violates Saito's criterion.
\end{pf}

\bibliographystyle{plain}\bibliography{algbib}

\end{document}